\documentclass[a4paper,draft, 10pt]{amsart}
\usepackage{amssymb, amsmath}
\usepackage{amscd} 
\usepackage{enumerate, multicol,amsfonts}
\usepackage{cite}

\theoremstyle{plain}
 \newtheorem{thm}{Theorem}[section]
  \newtheorem{definition}{Definition}[section]
 
 \newtheorem{lem}{Lemma}[section]
 
\theoremstyle{definition}
 
 \newtheorem{dfn}{Definition}[section]
\theoremstyle{remark}
 \newtheorem{rem}{Remark}[section]
 \numberwithin{equation}{section}

\renewcommand{\leq}{\leqslant}
\renewcommand{\geq}{\geqslant}
\renewcommand{\setminus}{\smallsetminus}
\title[Good Lambda Inequality]{ Good Lambda Inequalities for Non-doubling Measures in $\mathbb{R}^n$}



\author[Dr. M.Bhandari]{\bfseries Dr. Mukta Bhandari}

\address{
Department of Mathematics \\ 
Chowan University   \\ 
Murfreesboro, NC 27855\\
USA}
\email{bhandm@chowan.edu}

\begin{document}

\vspace{18mm} \setcounter{page}{1} \thispagestyle{empty}

\begin{abstract}
We establish a good lambda inequality
relating to the distribution function of Riesz potential and fractional
maximal function on $\left( \mathbb{R}^n, d\mu \right) $ where $\mu$ is a positive Radon
measure which doesn't necessarily satisfy a doubling condition. This
is extended to weights $w$ in $A_{\infty}(\mu)$ associated to the
measure $\mu$. We also derive potential inequalities as an
application.
\end{abstract}

\maketitle

\begin{center} \section{\large INTRODUCTION}
\end{center}
In this paper we will discuss the good-$\lambda$ inequality for
Riesz potentials associated to a Radon measure not necessarily
doubling but satisfies only a mild condition, which we call a growth
condition in $\mathbb{R}^n$.\\
\begin{definition}(~\cite{AH})  Let $\left(\mathbb{X}, d \right) $ be a metric measure space equipped with a metric $d$ and a Radon Measure $\mu$. The Riesz potential operator of order
	$\alpha $ on $\mathbb{X}$ is given by 
	
	\begin{equation} I_{\alpha}f(x) =
	\int_{\mathbb{X}} \frac{f(y)d\mu (y)}{d(x,y)^{N-\alpha}}, x\in
	\mathbb{X} \end{equation}
where $N$ is a fixed positive integer and $0<\alpha<N$.
\end{definition}
\begin{definition} (See ~\cite{CG})
	A Borel measure $\mu$ on a measure metric space $(\mathbb{X}, d)$ is
	said to satisfy the growth condition if  \begin{equation}
	\label{e:nondoubling}\mu \left( B(x,r)\right) \leq C r^N \end{equation} where
	the constant $C$ is independent of $x$ and $r$. This allows, in
	particular, non-doubling measures.
\end{definition}

\begin{definition}(See ~\cite{CG})
	A measure $\mu$ on a measure metric space $(\mathbb{X}, d)$ is said
	to satisfy the so-called ``doubling condition" if there exists a
	constant $C=C(\mu )\geq 1$, such that, for every ball $B(x,r)$ of
	center $x$ and radius $r$ 
	\begin{equation} \label{e:doubling}\mu
	\left(B(x,2r)\right) \leq C\mu\left(B(x,r)\right).\end{equation}
\end{definition}

\begin{definition}
	A pair of non-negative measurable functions $f$ and $g$ defined on
	$\mathbb{R}$ are said to satisfy a good-$\lambda$ inequality if 
	there exists a constants $K>1$, $0<\epsilon_{0}\leq 1$ such that for
	every $\lambda>0$
	
 \begin{equation} |\{x\in \mathbb{R}: f(x)>K\lambda,
	g(x)<\epsilon \lambda \}| \leq b(\epsilon)|\{x\in \mathbb{R}: f(x) >
	\lambda \}| 
\end{equation} 

where $b(\epsilon) \rightarrow 0$ as $\epsilon
	\rightarrow 0.$
\end{definition}

In their ground-breaking paper ~\cite{BG}, Burkholder and Gundy showed that the random variables $\mathbb{X}^{\star}$ and
$S(\mathbb{X})$ satisfy certain inequalities relating to their
distributions. These are now commonly called {\it good-$\lambda$
	inequalities}. In the harmonic function setting the first
good-$\lambda$ inequalities were also proved by Burkholder and Gundy
~\cite{BG1}. They were subsequently improved and refined by, among
others, Burkholder ~\cite{Bu3}, Dahlberg ~\cite{Da}, Fifferman,
Gundy, Silverstein and Stein ~\cite{FGSS}. The variations of
good-$\lambda$ inequalities  and its applications can be found in
the work of Benjamin Muckenhoupt and Richard L. Wheeden (see
\cite{MW}), S. D. Jaka (see \cite{J}), D. L. Burkholder (see
\cite{Bu1}), Rodrigo Ba\~{n}uelos (see \cite{B1}), and Richard F.
Bass (see \cite{B2}). A fair amount of deal about such inequalities
is found in the book by Rodrigo Ba\~{n}uelos and Charles N. Moore
(see \cite{BM}). Let us recall the classical good lambda
inequalities of Burkholder and Gundy ~\cite{BG} for continuous time
martingales.

\begin{thm} Let $\mathbb{X}_t$ be a continuous time martingale with
	maximal function $\mathbb{X}^{\star}$ and square function
	$S(\mathbb{X})$. Then for all $0<\epsilon<1$, $\delta>0$ and
	$\lambda >0$, 
\begin{displaymath} P\left\{\mathbb{X}^{\star}>\delta \lambda,
	S(\mathbb{X})\leq \epsilon \lambda \right\}\leq
	\frac{\epsilon^2}{(\delta-1)^2}P\{\mathbb{X}^{\star} >\lambda \} 
\end{displaymath}
	and 
\begin{displaymath} P\left\{ S(\mathbb{X})> \delta
	\lambda,\mathbb{X}^{\star}\leq \epsilon \lambda \right\}\leq
	\frac{\epsilon^2}{(\delta-1)^2}P\{S(\mathbb{X})
	>\lambda \} 
\end{displaymath}
\end{thm}

As they are expressed here, these are actually a refinement, due to
Burkholder (see \cite{Bu2}), of the inequalities of ~\cite{BG}. The
usefulness of such inequalities is already amply demonstrated by the
following lemma, which is but one of the many applications of these
type of inequalities. For this lemma we consider a  non-decreasing
function $\Phi$ defined on $[0, \infty]$ with $\Phi(0)=0$, $\Phi$ is
not identically $0$, and which satisfies the condition
$\Phi(2\lambda) \leq c\Phi(\lambda)$ for every $\lambda >0$, where
$c$ is a fixed constant. This lemma is from ~\cite{Bu2}.
\begin{lem}
	Suppose that $f$ and $g$ are nonnegative measurable functions on a
	measurable space $(\mathbb{Y}, \mathcal{A}, \mu)$, and $\delta >0$,
	$0<\epsilon <1$, and $0<\gamma <1$ are real numbers such that \begin{displaymath}
	\mu\{g>\delta \lambda, f\leq \epsilon \lambda \} \leq \gamma
	\mu\{g>\lambda\} \end{displaymath} for every $\lambda >0$. Let $\rho$ and $\nu$ be
	real numbers which satisfy \begin{displaymath} \Phi(\delta \lambda) \leq \rho
	\Phi(\lambda), \qquad \Phi(\epsilon^{-1}\lambda) \leq \nu
	\Phi(\lambda) \end{displaymath} for every $\lambda >0$. Finally, suppose $\rho
	\gamma < 1$ and $\int \Phi(\min \{1,g\})d\mu < \infty $. Then \begin{displaymath}
	\int_{\mathbb{Y}}\Phi(g)d\mu \leq \frac{\rho \nu}{1-\rho
		\gamma}\int_{\mathbb{Y}}\Phi (f) d\mu. \end{displaymath}
\end{lem}
We need the following definition for our purpose:

\begin{dfn}
	Let $(\mathbb{X},d)$ be a metric space and $\mu$ be any measure on
	$\mathbb{X}$. Let $f$ be a locally integrable function on
	$\mathbb{X}$. The maximal function of $f, \mathcal{M}(f)$, is
	defined by 
	\begin{displaymath} 
		\mathcal{M}(f)(x) =\sup_{r>0}\frac{1}{\mu
		\left(B(x,r)\right)}\int_{B(x,r)}|f(y)|d\mu(y)
	\end{displaymath} 
The operator
	$\mathcal{M}$ is called Hardy-Littlewood maximal operator. The
	fractional maximal function of $f$ is defined for $0<\alpha <N $ by
	\begin{displaymath} M_{\alpha} (x) = \sup_{r>0}
	\frac{1}{\mu(B(x,r))^{\frac{N-\alpha}{N}}}\int_{B(x,r)}|f(y)|d\mu(y).
	\end{displaymath}
\end{dfn}

Muckenhoupt and Wheeden in ~\cite{BM} proved the following
version of good $\lambda $- inequality in 1972.

\begin{thm}
	Let $\mu $ be a positive Radon measure in $\mathbb{R}^n $. Then there exists
	$a>1, b>0 $ such that for every $\lambda >0$ and for every $
	\epsilon, 0< \epsilon \leq 1$ \begin{displaymath} | \{x: I_{\alpha}\mu (x) >
	a\lambda \}| \leq b\epsilon ^{n/(n-\alpha)}| \{ x: I_{\alpha}\mu (x)
	> \lambda \}| + |\{ x: \mathcal {M}_{\alpha}\mu(x) > \epsilon
	\lambda \}| \end{displaymath} where \[I_{\alpha} \mu(x)= \int_{\mathbb{R}^n} \frac {d\mu
		(y)}{|x-y|^{n-\alpha}}.\] is the Riesz potential of the measure
	$\mu$.
\end{thm}

As a special case, the above theorem is still true if $d\mu(y)$ is
replaced by $f(y)dm(y)$ where $f$ is a measurable function and $m$
is the Lebesgue measure in $\mathbb{R}^n$. Note that Lebesgue measure is
doubling. In general, the above inequality is true for any positive
Radon measure which is doubling. The following theorem reveals this.

\begin{thm}\label{t:glifordoubling}
	Let $0<\alpha < N$, $\frac{1}{q}= 1-\frac{\alpha}{N}$, and $I_{\alpha}(x) =
	\int_{\mathbb{R}^n}\frac{f(y)d\mu(y)}{|x-y|^{N-\alpha}}$ where $\mu$ is a
	positive Radon measure satisfying the doubling condition
	\eqref{e:doubling} and $f$ is a measurable function on $\mathbb{R}^n$. Then
	there exist constants $a>1, b>0$ such that for every $\epsilon,
	0<\epsilon \leq 1$ and $\lambda
	>0$ we have \begin{displaymath} \mu \left(\{x:I_{\alpha}(x)
	>a\lambda\}\right) \leq b\epsilon^{N/N-\alpha}\mu
	\left(\{x:I_{\alpha} >\lambda\}\right)+\mu
	\left(\{x:M_{\alpha}f(x)>\epsilon \lambda\}\right).\end{displaymath}
\end{thm}
Our main objective in this paper is to establish the above
inequality for a Radon measure $\mu$ which is not necessarily
doubling but only satisfies a mild condition that we defined as
``{\it growth condition}'' in ~\eqref{e:nondoubling}.

\section{Useful Remarks and Results}
In this section, we will mention some useful remarks and results
in order to use them for our main result. We denote by $Q$ a cube in
$\mathbb{R}^n$ with sides parallel to the coordinate axes and by $\ell(Q)$
the side length of $Q$. Xavier Tolsa in ~\cite{To1} makes the
following remarks.

\begin {rem}If $\mu $ satisfies the growth condition ~\eqref{e:nondoubling} then there are lots
of big doubling cubes. Precisely speaking, given any point $x\in $
supp$(\mu )$ and $c>0$, there exists some $(\alpha, \beta
)$-doubling cube $Q$ (that is, $\mu (\alpha Q ) \leq \beta \mu (Q)$
where $\alpha >1$ and $\beta > \alpha ^n)$ centered at $x$ with
$\ell (Q) \geq c$. 
\end{rem}
\
\paragraph{} Note that if $\alpha, \beta $ are not specified then by a doubling cube we
will mean a $(2,\beta)$-doubling cube where $\beta > 2^n$.

\begin{rem}
There are always small doubling cubes for a Radon measure $\mu
$. That is, for $\mu$-a.e. $x\in \mathbb{R}^n$ there is a sequence of $(\alpha,
\beta)$-doubling cubes $\{Q_j\}_j^{\infty}$ centered at $x$ such
that $\ell(Q_j)\rightarrow 0$ as $j\rightarrow \infty.$
\end{rem}
Jos\'{e} Garc\'{\i}a-Cuerva and A. Eduardo Gatto proved the
following theorem in 2003 in ~\cite{CG}.

\begin{thm}\label{t:gatto} 
Let $(\mathbb{X},d,\mu)$ be a metric measure space where
the measure $\mu $ is Borel which satisfies the growth condition
~\eqref{e:nondoubling}. Then, for $1\leq p<\frac{n}{\alpha} $ and $
\frac{1}{q} = \frac{1}{p}-\frac{\alpha}{n}$, we have \[ \mu \left(\{
x\in \mathbb{X}:| I_{\alpha}f(x)|>\lambda \}\right)\leq
\left(\frac{C||f||_{L^p(\mu )}}{\lambda}\right)^q, \]that is ,
$I_{\alpha}$ is a bounded operator from $L^p(\mu)$ into the Lorentz
space $L^{q,\infty}(\mu )$.
\end{thm}

The above theorem establishes the weak type estimate for the Riesz
potentials in a metric space associated to Borel measure $\mu$ for
$1\leq p < \infty$.

\section{Main Results}
In this section, we will prove the good-$\lambda$ inequality for
Riesz potentials in $\mathbb{R}^n$ equipped with a measure not necessarily
doubling but satisfies only the mild condition called the ``growth
condition" ~\eqref{e:nondoubling}.
\begin{thm}
	Let $0<\alpha<N$. The Riesz potential operator of order
	$\alpha $ on $\mathbb{R}^n$ is given by \[ I_{\alpha}f(x) = \int_{\mathbb{R}^n}
	\frac{f(y)d\mu (y)}{d(x,y)^{N-\alpha}}, x\in \mathbb{R}^n \] for measurable
	function $f$ on $\mathbb{R}^n$. Let $\mu$ be a positive Radon measure on
	$\mathbb{R}^n$ satisfying the growth condition ~\eqref{e:nondoubling}. Then
	there exists constants $k\geq 1$ and  $C$ such that for every $\lambda >0$
	and $\epsilon , 0<\epsilon \leq 1,$\begin{equation} \label{e:gliforndoubling} \mu
	\left( \{x: I_{\alpha}f(x)>k\lambda, \mathcal{M}_{\alpha}f(x) \leq
	\epsilon \lambda \}\right) \leq C \epsilon^{\frac{N}{N-\alpha}}\mu
	\left(\{x: I_{\alpha }f(x) >\lambda \}\right).\end{equation}
\end{thm}

\begin{proof}
	
	Let $E_{\lambda } = \{ x\in \mathbb{R}^n : I_{\alpha} f(x) > \lambda \}.$
	Then $E_{\lambda }$ is open. So it has the following Whitney
	decomposition (see appendix J in ~\cite{Gr}).\\
	 There exists a
	countable family of dyadic cubes $\{Q_j\}_j^\infty$ such that
	
	\begin{enumerate}[(i)]
		\item $E_{\lambda} = \bigcup_j Q_j$ where $Q_j$'s have disjoint interiors.
		\item $diam(Q_j) \leq dist(Q_j, E_{\lambda} ^c) \leq 4 diam(Q_j)$, for
		every j. That is, 
	\begin{displaymath} \sqrt{n}\ell (Q) \leq dist(Q_j,
		E_{\lambda}^c) \leq 4\sqrt{n}\ell (Q). 
	\end{displaymath}
		\item If the boundaries of two cubes $Q_j$ and $Q_k$ touch, then
		
		\begin{displaymath} 
			\frac{1}{4} \leq \frac{\ell (Q_j)}{\ell (Q_k)}\leq 4 .
		\end{displaymath}
	
		\item For a given $Q_j$, there exists at most $12^n$ $Q_k$'s that
		touch it.
		\item For every $0<\delta < \frac{1}{4}, (1+\delta )Q \subseteq E_{\lambda
		},$ and  $\sum_{Q} \chi_{(1+\delta )Q}(x) \leq 12^n
		\chi_{E_{\lambda }} (x)$ for every $x\in E_{\lambda }$. 
		
		This implies that 
	\begin{displaymath} \sum_{Q} \mu \left( (1+\delta )Q \right)(x)
		\leq 12^n \mu (E_{\lambda }).
	\end{displaymath}
	\end{enumerate}

	Fix $x\in E_{\lambda}$ and a $Q$ in $\{Q_j\}_j^\infty$ with $x\in
	Q$.  Let $k> 1$ and consider the set $\{x\in Q: I_{\alpha}(x) > k \lambda
	\}$. Suppose $Q\bigcap \{\mathcal{M}_{\alpha}f(x) \leq \epsilon
	\lambda \}\neq \phi$.  Let $Q_x$ be the cube center at $x$ and
	$\ell(Q_x)=\delta \ell (Q)$. We may fix $\delta = \frac{1}{8}.$ Then
	$Q_x \subseteq (1+\delta )Q$ and $16Q_x \supseteq 2Q$. Consider
	$\frac{1}{2^j}Q_x, j=1,2,3,\ldots$. Take the first which is doubling
	and denote it by $\widehat{Q}_x:= \frac{1}{2^{m+1}}Q_x$ where $m$
	depends on $x$ and $\widehat{Q}_x \subseteq \frac{1}{2}Q_x$. Thus,

	\begin{enumerate}[(a)]
		\item For every $x\in E_{\lambda}\bigcap \{\mathcal{M}_{\alpha}f(x)\leq \epsilon \lambda \}$
		there exists a doubling cube $\widehat{Q}_x$ such that $\mu
		(2\widehat{Q}_x ) \leq \beta \mu (\widehat{Q}_x)$ where $\beta >
		2^n$. We may assume $\beta = 2^{n+\epsilon}.$ \vspace{.2in}
		\item $\mu \left( 2.\frac{1}{2^{j+1}}Q_x\right) \geq \beta \mu \left(
		\frac{1}{2^{j+1}}Q_x\right)$ for every $j<m.$  \vspace {.2in}
	\end{enumerate}

	Therefore, 
	\begin{equation*} 
		\{x\in Q: I_{\alpha}f(x)>k\lambda,
	\mathcal{M}_{\alpha}f(x)<\epsilon \lambda \} \subseteq \bigcup_{x\in
		Q}\widehat{Q}_x.
	 \end{equation*} 
 
 The Besicovitch Covering Lemma implies that
	there exists $\widehat{Q}_{x_j}, j=1,2,3,\ldots$ such that 
	
	\begin{displaymath} \{
	x\in Q:I_{\alpha}f(x) >k \lambda, \mathcal{M}_{\alpha}f(x) <
	\epsilon \lambda \} \subseteq \bigcup_{j}\widehat{Q}_{x_j} 
\end{displaymath}
 where
	\begin{displaymath} 
		\sum_{j}\chi _{\widehat{Q}_{x_j}} \leq 4^n
	\chi_{(1+\delta)Q}.
\end{displaymath} That is, 
\begin{displaymath} \mu \left(
	\bigcup_{j}\widehat{Q}_{x_j}\right) \leq \mu \left(
	(1+\delta)Q\right).
 \end{displaymath}
	Let $f=f_1 +f_2$ where $f_1=\chi_{2Q}f$ and $f_2=f-f_1.$
	Then,
	\begin{align*} \{x \in Q:I_{\alpha}f_1(x)>k\lambda,&
	\mathcal{M}_{\alpha}f(x)<\epsilon \lambda\} \\ 
	&\subseteq
	\bigcup_{j}\{x\in \widehat{Q}_{x_j}:I_{\alpha}f_1(x)>k\lambda,
	\mathcal{M}_{\alpha}f(x)<\epsilon \lambda\}.
\end{align*}

 Fix a doubling cube
	$\widehat{Q}_{x_j}$. Write $f_1 =\chi_{2\widehat{Q}_{x_j}}f_1 +
	\chi_{(2\widehat{Q}_{x_j})^c}f_1=f_{11}+f_{12}.$ Then, by using the
	weak type inequality (see theorem \eqref{t:gatto}, ~\cite{CG}), we
	obtain
	
		\begin{align*}
		\mu (\{ x\in 2\widehat {Q}_{x_j}: & I_\alpha f_{11} >  k\lambda \})  \\ 
		& \leq \left( \frac{1}{k\lambda}\|f_{11}\|_1 \right)^{N/N-\alpha}\\
		&= \left( \frac{1}{k\lambda}\int_{2\widehat{Q}_{x_j}}|f_1(x)|d\mu(x)\right) ^{N/N-\alpha}\\
		&=\left(\frac{1}{k\lambda}\mu(2\widehat{Q}_{x_{j}})^{\frac{N-\alpha}{N}}\frac{1}{\mu(2\widehat{Q}_{x_j})^{\frac{N-\alpha}{N}}}
		\int_{2\widehat{Q}_{x_j}}|f_1(x)|d\mu(x)\right)^{N/N-\alpha}\\
		&= \left(\frac{1}{k\lambda}\right)^{\frac{N}{N-\alpha}}\mu(2\widehat{Q}_{x_j}) M_\alpha f_1(x_j)^{N/N-\alpha}\\
		&\leq C\epsilon^{N/N-\alpha}\mu(\widehat{Q}_{x_j}).
		\end{align*}

 Note that $\mu(Q) \leq c\ell(Q)^N$ for any cube $Q$. From this, it follows that
	$$\left(\frac{c}{\mu(Q)}\right)^{\frac{N-\alpha}{N}} \geq
	\frac{1}{\ell(Q)^{N-\alpha }}.$$  Also,
	 \begin{equation*} \mu(Q_{x_j}) \geq \beta
	\mu\left(\frac{1}{2}Q_{x_j}\right) \geq \ldots > \beta^{i}
	\mu\left(\frac{1}{2^i}Q_{x_j}\right)
\end{equation*}
 where $i<m$. This follows
	by the choice of our $Q_{x_j}$'s. So, if $x\in \widehat {Q}_{x_j}$
	
	\begin{align*}
		I_\alpha f_{12}(x) &= \int_{(2\widehat{Q}_{x_j})^c} \frac{f_1(y)}{d(x,y)^{N-\alpha}}d\mu(y)\\
		&\leq \sum_{k=0}^{m-1}\int_{\frac{1}{2^k}{Q_{x_j}}\setminus\frac{1}{2^{k+1}}{Q_{x_j}}}
		\frac{|f_1(y)|}{d(x,y)^{N-\alpha}}d\mu(y)+
		\int_{16Q_{x_j}\setminus 2Q}\frac{|f_1|(y)d\mu(y)}{d(x,y)^{N-\alpha}}\\
		&\leq \sum_{k=0}^{m-1}\frac{c}{\ell(\frac{1}{2^k}Q_{x_j})^{N-\alpha}}\int_{\frac{1}{2^k}{Q_{x_j}}
			\setminus\frac{1}{2^{k+1}}{Q_{x_j}}}|f_1(y)|d\mu(y) \\
		&+\frac{C}{\ell (16Q_{x_j})^{N-\alpha}}\int_{16Q_{x_j}\setminus 2Q}f_1(y)d\mu(y)\\
		&\leq \sum_{k=0}^{m-1}\frac{2^{k(N-\alpha)}}{\ell (Q_{x_j})^{N-\alpha}}\int_{\frac{1}{2^k}
			{Q_{x_j}}\setminus\frac{1}{2^{k+1}}{Q_{x_j}}}|f_1(y)|d\mu(y) \\ &+
		\frac{C}{\mu (16Q_{x_j})^{(N-\alpha)/N}}\int_{16Q_{x_j}\setminus 2Q}f_1(y)d\mu(y)\\
		&\leq \sum_{k=0}^{m-1}C2^{k(N-\alpha)}\frac{1}{\mu(Q_{x_J})^{\frac{N-\alpha}{N}}}
		\int_{\frac{1}{2^k}{Q_{x_j}}\setminus\frac{1}{2^{k+1}}{Q_{x_j}}}|f_1(y)|d\mu(y)+CM_{\alpha}f_1(x_j)\\
		&\leq \sum_{k=0}^{m-1}C2^{k(N-\alpha)}\frac{1}{\beta^k\mu(\frac{1}{2^k}Q_{x_j})^{\frac{N-\alpha}{N}}}
		\int_{\frac{1}{2^k}{Q_{x_j}}\setminus\frac{1}{2^{k+1}}{Q_{x_j}}}|f_1(y)|d\mu(y)+ C\lambda\\
		&\leq \sum_{k=0}^{m-1}C\frac{2^{k(N-\alpha)}}{\beta^{\frac{k}{N}(N-\alpha)}}M_{\alpha}f_1(x_j)+C\lambda\\
		&= C\sum_{k=0}^{m-1}\left(\frac{2^{N-\alpha}}{2^{N-\alpha}2^{\epsilon(\frac{N-\alpha}{N})}}\right)^{k}
		M_{\alpha}f_1(x_j)+C\lambda\\
		&= C\sum_{k=0}^{m-1}2^{-\epsilon(\frac{N-\alpha}{N})k}M_\alpha f_1(x_j)+C\lambda\\
		&=C_{N,\alpha}M_{\alpha}f_1(x_j) \leq C_{N,\alpha}=C\lambda.
	\end{align*}

	Note that the constant C in different occurrences above are not
	necessarily the same. Therefore 
	
	\begin{displaymath} 
		I_{\alpha}f_{12}(x) \leq C\lambda
\end{displaymath} 

for every $x\in \widehat{Q}_{x_j}$. We now estimate
	$I_{\alpha}f_2$. Fix a point $x\in Q$. Consider the ball $B=B(x,
	6diam(Q)) \supseteq 2Q$. Let $x_0\in B\bigcap{E_{\lambda}}^c$ such
	that 
	
	\begin{displaymath} diam(Q) \leq dist(Q,x_0)\leq 4diam(Q).
	\end{displaymath} 

Then for any
	$y\in (2Q)^c$
	
	\begin{align*}
		d(x_0,y) &\leq d(x,x_0)+d(x,y)\\ &\leq C diam(Q)+d(x,y) \\ &\leq C
		dist(Q,E_{\lambda}^c) +d(x,y) \\ &\leq Cd(x,y)+ d(x,y).
	\end{align*}

	Thus 
	\begin{displaymath} d(x_0,y)\leq Cd(x,y)
	\end{displaymath} 
for $x_0\in B\bigcap E_{\lambda}^c$
	and $x\in Q.$  Therefore, 
	
	\begin{displaymath}
	I_{\alpha}f_2(x)=\int_{(2Q)^c}\frac{f(y)d\mu(y)}{d(x,y)^{N-\alpha}}\leq
	CI_{\alpha}f(x_0)<C\lambda.
\end{displaymath}

 Finally, summing over all Whitney cubes yields the inequality. Indeed, 
	
	\begin{align*}
		\mu &\left( \{ x:I_\alpha f(x) >k\lambda, M_\alpha f(x) <\epsilon \lambda \}\right)\\
		& =\sum_{Q} \mu \left( \{ x\in Q:I_\alpha f(x) >k\lambda, M_\alpha f(x) <\epsilon \lambda \}\right)\\
		&\leq \sum_{Q} \mu\left(\{ x\in Q : I_\alpha
		f_1  >
		(k-C)\lambda, M_\alpha f < \epsilon \lambda \}\right)\\
		&= \sum_{Q} \sum _{j} \mu\left(\{ x\in{\widehat{Q}_{x_j}} :
		I_\alpha f_1  >
		(k-C)\lambda, M_\alpha f < \epsilon \lambda \}\right)\\
		&<\sum_{Q}\sum_{j}\mu\left(\{ x\in{\widehat{Q}_{x_j}} :
		I_\alpha f_{11}  >
		(k-2C)\lambda, M_\alpha f < \epsilon \lambda \}\right)\\
		& \leq \sum_{Q}\sum_{j}C\epsilon^{N/N-\alpha}\mu(\widehat{Q}_{x_j})\\
		&\leq \sum_{Q} C\epsilon^{N/N-\alpha}4^n \mu((1+\delta)Q)\\
		&= 4^n 12^n C \epsilon ^{N/N-\alpha} \mu(\{I_\alpha f > \lambda
		\}).
	\end{align*}

	Thus, finally we have 
	
	\begin{align*}
		\mu\left(\{ x: I_\alpha f > k\lambda \}\right) & \leq \mu \left(
		\{ I_\alpha f > k\lambda , M_\alpha f \leq \epsilon \lambda
		\}\right)
		+ \mu \left(\{ M_\alpha f > \epsilon \lambda \}\right) \\
		&\leq C\epsilon ^{N/N-\alpha} \mu \left( \{ I_\alpha f > \lambda
		\}\right) + \mu \left( \{ M_\alpha f > \epsilon \lambda
		\}\right)
	\end{align*}

	where the constant $C=C(n,N,\alpha )$. This completes the proof of the theorem. 
\end{proof}

\section{\bf Good-$\lambda$ inequality for weights}

Next, we extend this inequality
for the weights
$w\in A_{p}(\mu)$ where $\mu$ satisfies the growth condition
~\eqref{e:nondoubling}. That is for any measurable set $E$,\begin{displaymath}
w(E)=\int_{E}w(x)d\mu(x),\end{displaymath} where $\mu$ is a positive Radon
measure which satisfies the growth condition \eqref{e:nondoubling}.
Details about this type of weights can be obtained from the paper
~\cite{OP} by Joan Orobitg and Carlos P\'{e}rez.
$A_{\infty}(\mu)$ is defined as $A_{\infty}(\mu) =
\bigcup_{p>1}A_p(\mu)$ in a classical way.

\begin{thm}
	Let $w\in A_{\infty}(\mu)$, $\alpha >0$, and $0<\lambda <\infty$.
	Then there exists a positive constant $a\geq 1$ for which, for every
	$\eta >0$, there is an $\epsilon, 0<\epsilon \leq 1$, such that the
	inequality \begin{displaymath} \begin{split} \label{e:gliforweight} &w\left(\{x\in
		\mathbb{R}^n : I_{\alpha}f(x)
		> a\lambda\}\right)\\ &\leq \eta w \left(\{x\in \mathbb{R}^n :I_{\alpha}f(x) >
		\lambda \}\right) + w \left(\{x\in \mathbb{R}^n: M_{\alpha}f(x)>\epsilon
		\lambda \}\right) \end{split} \end{displaymath}  holds for every $\lambda >0.$
\end{thm}
\begin{proof}
	Let $E_{\lambda} = \{x\in \mathbb{R}^n: I_{\alpha}(x) >\lambda \},
	\lambda >0$. Then $E_{\lambda} $ is open because $I_{\alpha} $ is lower
	semi-continuous. Then there exists a family of dyadic cubes
	$\{Q_j\}$, called Whitney cubes, such that $E_{\lambda} =
	\bigcup_{j}Q_j $ and 
	\begin{displaymath} diam(Q_j) \leq dist(Q_j, E_{\lambda}) \leq
	4 diam(Q_j).\end{displaymath} 

Because $w\in A_{\infty}(\mu ) $, it follows that,
	for every $\eta >0$, there is a $\delta $ such that if $Q$ is a cube
	and $E$ is a measurable subset of $Q$ then there is a constant $C_0$
	such that 
	\begin{displaymath} \frac{w(E)}{w(Q)}\leq C_0 \left(\frac{\mu (E)}{\mu
		(Q)}\right)^{\delta}. \end{displaymath} (see ~\cite{OP}).\vspace{.25in} That
	is for every $Q\in \{Q_{j}\}$,
	
	\begin{align*} 
	\frac {w \left( \{x\in Q:  I_{\alpha}f(x) >k \lambda, M_{\alpha}f(x) <\epsilon
		\lambda \}\right)}{w(Q)} & \leq  C_0 \left( \frac{\mu \left( \{x\in Q:
		I_{\alpha}f(x)>k \lambda, M_{\alpha}f(x) <\epsilon \lambda
		\}\right)}{\mu(Q)}\right)^{\delta}\\ \leq & C_0 \left( C_1\epsilon ^{N/N-\alpha}\right)^{\delta}.
\end{align*}

	This implies that 
	\begin{displaymath} 
	w \left( \{x\in Q:  I_{\alpha}f(x) >k \lambda, M_{\alpha}f(x)
	<\epsilon \lambda \}\right)  \leq C\epsilon ^{\left(N/N-\alpha
		\right)\delta}w(Q).
	\end{displaymath} 

It then follows that 

\begin{displaymath} w\left( \{ x\in Q:
	I_{\alpha}(x) > k\lambda\}\right)\leq C\epsilon ^{\left(N/N-\alpha
		\right)\delta}w(Q) + w\left(\{x\in Q: M_{\alpha}(x) > \epsilon \lambda
	\}\right).
\end{displaymath}

	The theorem now follows from summing over all $Q\in
	\{Q_j\}.$
\end{proof}
\paragraph{} We observe that the above good lambda inequality
is true for the weights $w\in A_{\infty}$ associated to a measure
$\mu$ which is doubling. This is because the good-$\lambda$
inequality holds for doubling measure as well (see theorem
~\eqref{t:glifordoubling}).

\section{Applications} 

In this section we will provide some
applications of our main results. Basically, we will derive the norm
inequalities for fractional integrals and maximal functions. Note
that the associated measure $\mu$ that follows in this section
satisfies either the growth condition ~\eqref{e:nondoubling} or the
doubling condition ~\eqref{e:doubling}. This is because our main
results are true in either case. Here follows the two theorems as an application of our main results.

\begin{thm}
	Let $1<p<\infty$ and $0<\alpha < N$. Then there is a constant $A$
	such that for any measurable function $f$ on $\mathbb{R}^n$ and a Radon
	measure $\mu$ satisfying either the growth condition
	~\eqref{e:nondoubling} or the doubling condition
	~\eqref{e:doubling}, we have 
	
	\begin{displaymath} \|I_{\alpha}f\|_{p}\leq
	A\|M_{\alpha}f\|_{p}. 
\end{displaymath}

\end{thm}

\begin{proof}
	 We assume that $\mu$ has compact support. The right
	hand inequality is a consequence of "good lambda inequality"
	\eqref{e:gliforndoubling}. Multiplying the good lambda inequality
	\eqref{e:gliforndoubling} by $\lambda ^{p-1}$, we obtain     for any
	positive $R$,
	
	\begin{align*}
		\int_{0}^R \mu (\{x:&I_{\alpha}f(x) >k \lambda\})\lambda^{p-1}d\lambda \\
		& \leq b\epsilon^{N/N-\alpha} \int_0^R \mu (\{ x: I_{\alpha}f(x) > \lambda \})\lambda^{p-1}d\lambda \\ 
		&+ \int_0^R \mu (\{ x: M_{\alpha}f(x) >\epsilon \lambda \})\lambda^{p-1} d\lambda.
	 \end{align*}
 
	After changing variables we obtain,
	
	\begin{align*}
		k^{-p} &\int_0^{kR} \mu(\{I_{\alpha}f(x) > \lambda \})\lambda^{p-1}d\lambda \\
		& \leq \frac{1}{2} k^{-p} \int_0^R \mu(\{ I_{\alpha}f(x)> \lambda \}) \lambda ^{p-1}d\mu +
		\epsilon^{-p} \int_0^{kR} \mu ( \{ M_{\alpha}f(x) > \lambda \})\lambda ^{p-1}d\lambda.
	\end{align*}

	That is, 
	\begin{displaymath} k^{-p}\int_0^{kR} \mu(\{I_{\alpha}f(x) > \lambda
	\})\lambda^{p-1}d\lambda \leq 2\epsilon^{-p} \int_0^{kR} \mu ( \{
	M_{\alpha}f(x) > \lambda \})\lambda ^{p-1}d\lambda.
\end{displaymath}

 Letting $ R\rightarrow \infty $ and using the definition $$ \int_{X}|f|^pd\mu =
	p\int_0^\infty t^{p-1}\mu(\{ x: f(x) >t\}) dt $$ we obtain,
	
	 \begin{displaymath}
	k^{-p}\int_X |I_{\alpha}f(x)|^p d\mu (x) \leq 2\epsilon^{-p} \int_X
	|M_{\alpha}f(x)|^{p}d\mu 
\end{displaymath}

 This yields the right hand inequality. If $\mu$
	doesn't have compact support, we let $\mu_n$ be the restriction of
	$\mu$ to the ball $B(x_0,n)$ for $ n=1,2,3,\ldots$ where $x_0$ is
	some point in $X$. This yields $\|I_{\alpha}f\|_{L^p(X,\mu_n)} \leq A
	\|M_{\alpha}f\|_{L^p(X,\mu_n)} $ for all $n$, where $A$ is
	independent of $n$. Then the  theorem follows  by Monotone Convergence
	Theorem. 
\end{proof}

\begin{thm}
	Let $\mu$ be a positive radon measure satisfying either the doubling
	condition \eqref{e:doubling} or the growth condition
	\eqref{e:nondoubling}, $w\in A_{\infty}(\mu), 0<\alpha < N$, and let
	$0<p<\infty.$ Then there exists a positive constant $C$, that only
	depends on $N, p,$ and $A_{\infty}(\mu)$ constants of $w$, such that
	
	\begin{displaymath} \int_{\mathbb{R}^N}|I_{\alpha} |^pw d\mu \leq C\int_{\mathbb{R}^N}\left(M_{\alpha} \right)^pwd\mu 
	\end{displaymath} 

for every measurable function $f$.
\end{thm}

\begin{proof}
	Without loss of generality, we may assume that $f\geq 0$. We
	multiply the  inequality \eqref{e:gliforweight} by $\lambda^{p-1}$
	and integrate from 0 to $R (R>0)$ with respect to $\lambda$  to
	obtain

	\begin{align*} \int_0^R  w \left( \{x\in \mathbb{R}^N : I_{\alpha}f(x)
	> a\lambda\}\right) & p\lambda^{p-1} d\lambda  \\ \leq & \eta \int_0^R w \left(\{x\in \mathbb{R}^N :I_{\alpha}f(x) >
	\lambda \}\right)p\lambda^{p-1} d\lambda \\ & + \int_0^R w \left(\{x\in
	\mathbb{R}^N: M_{\alpha}f(x)>\epsilon \lambda \}\right)p\lambda^{p-1}
	d\lambda .
\end{align*}

	Applying change of variable and $a>1$ yields
	
\begin{align*}  a^{-p}  \int_0^{aR}   w ( \{x\in
	\mathbb{R}^N : I_{\alpha}f(x) 
	> & a \lambda\})  \lambda^{p-1} d\lambda  \\& \leq  \eta \int_0^{aR} w \left(\{x\in \mathbb{R}^N:I_{\alpha}f(x) >
	\lambda \}\right)\lambda^{p-1} d\lambda \\ & + \epsilon^{-p}\int_0^{\epsilon R} w \left(\{x\in
	\mathbb{R}^N: M_{\alpha}f(x)>\epsilon \lambda \}\right)\lambda^{p-1}
	d\lambda .
\end{align*}

	 Now we choose $ \eta \leq \frac{1}{2}a^{-p}$. This yields
	 
	\begin{align}  a^{-p}  \int_0^{aR}   w ( \{x\in
	\mathbb{R}^N : I_{\alpha}f(x)
		& > a\lambda\}) \lambda^{p-1} d\lambda \notag  \\ \leq &  2\epsilon^{-p}\int_0^{\epsilon R} w \left(\{x\in
		\mathbb{R}^N: M_{\alpha}f(x)>\epsilon \lambda \}\right)\lambda^{p-1}
		d\lambda .
	\end{align}

 Now let 
 
 \begin{displaymath} \chi(x,\lambda) =
	\begin{cases}
		1,  & \text{if $I_{\alpha}f(x) >\lambda >0$;}\\
		0, &\text{otherwise}
	\end{cases}
	\end{displaymath}
 Then,

 \begin{align*}  
	\int_0^{aR}w\left( \{x\in \mathbb{R}^N: i_{\alpha}f(x)
	>\lambda \}\right) & \lambda ^{p-1}d\lambda \\
	&= \int_0^{aR}\left(\int_{ \{x\in \mathbb{R}^N: I_{\alpha}f(x) > \lambda \}
	} w(x)d\mu(x)\right) \lambda^{p-1} d\lambda. \\
	&= \int_0^{aR}\left(\int_{\mathbb{R}^N}\chi(x,\lambda) w(x)d\mu(x)\right)
	\lambda^{p-1} d\lambda.\\
	&= \int_{\mathbb{R}^N} w(x)\int_0^{\min \{ aR,I_{\alpha}f(x) \}}\lambda
	^{p-1}d\lambda d\mu\\
	&= \int_{\mathbb{R}^N} \left( \min \{ aR, I_{\alpha}f(x) \}\right)^p w(x)d\mu (x).
\end{align*}

	Similarly, 
	
	\begin{displaymath}  \int_0^{\epsilon R}w\left( \{x\in \mathbb{R}^N: M_{\alpha}f(x)
	>\lambda \}\right)\lambda ^{p-1}d\lambda = \int_{\mathbb{R}^N} \left( \min \{ \epsilon R,  M_{\alpha}f(x) \}\right)^p w(x)d\mu
	(x). \end{displaymath}

 Therefore, 
 \begin{displaymath}  
 	\int_0^{aR}w\left( \{x\in \mathbb{R}^N:  I_{\alpha}f(x)
	>\lambda \}\right)\lambda ^{p-1}d\lambda \leq C \int_{\mathbb{R}^N} \left( \min \{ aR,  I_{\alpha}f(x) \}\right)^p w(x)d\mu
	(x).\end{displaymath}

 Using Fatou's lemma,

	\begin{align*} 
		\int_{\mathbb{R}^N} \left( I_{\alpha}f(x) \right)^p w(x)d\mu(x)&= \int_{\mathbb{R}^N} \liminf_{R\rightarrow
			\infty} \left( \min \{ aR, I_{\alpha}f(x) \}\right)^p w(x) d\mu(x) \\
		& \leq  \liminf_{R\rightarrow
			\infty}\int_{\mathbb{R}^N}  \left( \min \{ aR, I_{\alpha}f(x) \}\right)^p w(x) d\mu(x) \\
		&\leq  C \liminf_{R\rightarrow \infty} \int_{\mathbb{R}^N}  \left( \min \{ \epsilon R,  M_{\alpha}f(x) \}\right)^p w(x) d\mu(x) \\
		& \leq  \int_{\mathbb{R}^N} \left( M_{\alpha}f(x) \right)^p w(x)
		d\mu(x).
	\end{align*}

\end{proof}

\section{ Future Motivation} Muckenhoupt and Wheeden in ~\cite{MW}
have proved a weighted version of Sobolev imbedding theorem
as an application of good-$\lambda$
inequality associated to a Lebesgue measure. Following their
footprints, it is expected to obtain the Sobolev imbedding theorem
in more general context associated to a Radon measure which
satisfies either the growth condition \eqref{e:nondoubling} or the
doubling condition \eqref{e:doubling}. This will open a wide range
of spectrum for the future research on this line.

\bibliographystyle{IEEEbib}
\bibliography{egbib}

\end{document}